\documentclass{article}
\usepackage[utf8]{inputenc}
\usepackage{amsfonts}
\usepackage{amsmath, bm}
\usepackage{graphicx}
\usepackage{xcolor}
\usepackage{tikz}
\usepackage{hyperref}

 \newtheorem{theorem}{Theorem}[section]
 \newtheorem{lemma}{Lemma}[section]
\newtheorem{corollary}{Corollary}[section]
 \newtheorem{definition}{Definition}[section]
 \newtheorem{remark}{Remark}[section]

\newcommand{\xx}{\mathbf{x}}
\newcommand{\uu}{\mathbf{u}}
\newcommand{\ff}{\mathbf{f}}

\newcommand{\ZZ}{{\mathbb Z}}
\newcommand{\NN}{{\mathbb N}}

\title{Modeling excitable cells with the EMI equations: spectral analysis and iterative solution strategy}
\author{Pietro Benedusi, Paola Ferrari, Marie E. Rognes, \\ Stefano Serra-Capizzano}
\date{}

\begin{document}

\maketitle

\begin{abstract}
In this work, we are interested in solving large linear systems stemming from the Extra-Membrane-Intra (EMI) model, which is employed for simulating excitable tissues at a cellular scale. After setting the related systems of partial differential equations (PDEs) equipped with proper boundary conditions, we provide numerical approximation schemes for the EMI PDEs and focus on the resulting large linear systems. We first give a relatively complete spectral analysis using tools from the theory of Generalized Locally Toeplitz matrix sequences. The obtained spectral information is used for designing appropriate (preconditioned) Krylov solvers.  We show, through numerical experiments, that the presented solution strategy is robust w.r.t. problem and discretization parameters, efficient and scalable. 
\end{abstract}

\section{Introduction}
The EMI (Extra, Intra, Membrane) model, also known as the cell-by-cell model, is a numerical building block in computational electrophysiology, employed to simulate excitable tissues at a cellular scale. With respect to homogenized models, such as the well-known mono and bidomain equations, the EMI model explicitly resolves cells morphologies, enabling detailed biological simulations. For example, inhomogeneities in ionic channels along the membrane, as observed in myelinated neuronal axons \cite{black1990ion}, can be described within the EMI model. 

The main areas of application for the EMI model are computational cardiology and neuroscience, where spreading excitation by electrical signalling plays a crucial role \cite{buccino2021improving, ellingsrud2021cell,jaeger2021efficient,ellingsrud2020finite,jaeger2021millimeters,de2023boundary,tveito2017cell,mori2009numerical}.
We refer to \cite{EMI_book} for an exhaustive description of the EMI model, its derivation, and relevant applications.

This work considers a Galerkin-type approximation of the underlying system of partial differential equations (PDEs). A finite element discretization leads to block-structured linear systems of possibly large dimensions. The use of ad hoc tools developed in the theory of Generalized Locally Toeplitz matrix sequences allows us to give a quite complete picture of the spectral features (the eigenvalue distribution) of the resulting matrices and matrix sequences.

The latter information is crucial to understand the algebraic properties of such discrete objects and to design iterative solution strategies, enabling large-scale models, as in \cite{benedusi2018space, benedusi2021fast}. More precisely, spectral analysis is employed for designing tailored preconditioning strategies for Krylov methods. The resulting efficient and robust solution strategies are theoretically studied and numerically tested.

The current paper is organized as follows. In Section \ref{sec:EMI-pb} the continuous problem is introduced together with possible generalizations. 
Section \ref{sec:approx} considers a basic Galerkin strategy for approximating the examined problem. The spectral analysis is given in Section \ref{sec:spectral} in terms of distribution results and degenerating eigenspaces. Specifically, Subsection \ref{ssec:GLTbackground} lays out the foundational theories and concepts necessary for understanding the distribution of the entire system (presented in a general form in Subsection \ref{ssec:symbol_analysis}) and the specific stiffness matrices and matrix sequences (discussed in Subsection \ref{ssec:stiffness_matrices}). These findings are then employed in Subsection \ref{ssec:preconditioner} for proposing a preconditioning strategy. Moreover, in Section \ref{sec:num} we report and critically discuss the numerical experiments and Section \ref{sec:final} contains conclusions and a list of relevant open problems.

\section{The EMI problem}\label{sec:EMI-pb}

We introduce the partial differential equations characterizing the EMI model. When comparing to the homogenized models, we observe that the novelty of the EMI approach consists in the fact that the cellular membrane $\Gamma$ is explicitly represented, as well as the intra- and extra-cellular quantities, denoted with subscript $i$ and $e$ respectively. 

More in detail, given a domain $\Omega=\Omega_i\cup\Omega_e\cup\Gamma\subset\mathbb{R}^d$, typically with $d\in\{2,3\}$, $\partial\Omega=\partial\Omega_e/\Gamma$, and $\partial\Omega_i=\Gamma$, we consider the following time-discrete problem for the intra- and extra-cellular potentials $u_i,u_e$, and for the membrane current $I_m$:
\begin{align}\label{eq::EMI}
    & -\nabla \cdot (\sigma_e\nabla u_e(\xx))=0 && \quad \mathrm{for} \quad \xx\in\Omega_e, \\ 
    & -\nabla \cdot (\sigma_i\nabla u_i(\xx))=0 && \quad \mathrm{for} \quad \xx\in\Omega_i, \label{eq::EMI_2}\\
    & \sigma_e\nabla u_e(\xx)\cdot\mathbf{n}_e = -\sigma_i\nabla u_i(\xx)\cdot\mathbf{n}_i\equiv I_m(\xx) && \quad \mathrm{for} \quad \xx\in\Gamma, \label{eq::EMI_3} \\
    & u_i(\xx)-u_e(\xx) = v(\xx) && \quad \mathrm{for} \quad \xx\in\Gamma, \label{eq::EMI_4} \\ 
    & v(\xx)-\tau I_m(\xx)=f(\xx) && \quad \mathrm{for} \quad \xx\in\Gamma,\label{eq::EMI_5}
\end{align}
with $\sigma_e,\sigma_i,\tau\in\mathbb{R}_+$, $\mathbf{n}_i$ (resp. $\mathbf{n}_e$) is the outer normal on $\partial\Omega_i$ (resp. $\partial\Omega_e$) and $f\in L^2(\Gamma)$ is known. We can close the EMI problem with homogeneous boundary conditions:
\begin{align}\label{bc1}
     u_e(\xx) = 0 & \quad \mathrm{for} \quad \xx\in\partial\Omega_D,\\ \label{bc2}
     \sigma_e \nabla u_e(\xx) \cdot \mathbf{n}_e = 0 & \quad \mathrm{for} \quad \xx\in\partial\Omega_N,
\end{align}
with $\partial\Omega=\partial\Omega_D\cup\partial\Omega_N$. In the case of a pure Neumann problem, with $\partial\Omega=\partial\Omega_N$, uniqueness must be enforced through an additional constraint (e.g. a condition on the integral of $u_e$). We are imposing homogeneous boundary conditions for simplicity; the inhomogeneous case reduces to the homogeneous case by considering the lifting of the boundary data.

Essentially, the EMI problem consists of two homogeneous Poisson problems coupled at the interface $\Gamma$ through a Robin-type condition \eqref{eq::EMI_3}-\eqref{eq::EMI_5}, depending on the current $I_m$. The EMI problem can be considered a \textit{mixed-dimensional} problem since the unknowns of interest are defined on sets of different dimensionality, $d$ for $u_i$ and $u_e$ and $d-1$ for $I_m$.

It is worth noticing that possible dynamics originate in the membrane via the source term $f$. In particular, eq. \eqref{eq::EMI_5} is obtained by discretizing point-wise the capacitor current-voltage relation, a time-dependent ODE with $t\in(0,T]$ given a final time $T>0$:
\begin{align*}
   \frac{\partial v(t)}{\partial t} & = C^{-1}_m(I_m(t) - I_{\mathrm{ion}}(v,t)),\\
   v(0) &= v_0, \notag
\end{align*}
where the positive constant  $C_m$ is a capacitance and the ionic current $I_{\mathrm{ion}}$ is a reaction term. An implicit (resp. explicit) integration of $I_m$ (resp. $I_{\mathrm{ion}}$), with time step $\Delta t>0$, results in eq. \eqref{eq::EMI_5} with
\begin{equation}\label{eq::f}
f(\xx) = v_0(\xx) - \tau I_{ion}(v_0(\xx)),    
\end{equation}
and $\tau = C_m^{-1}\Delta t$. 

The model in eq. \eqref{eq::EMI}-\eqref{eq::EMI_5} can describe a single cell or multiple disjoint cells, i.e. with $\Omega_i=\bigcup_{j=1}^{N_\text{cell}}\Omega_i^j$, in an extracellular media, cf. Figure~\ref{fig:EMI}. 
Extending model \eqref{eq::EMI}-\eqref{eq::EMI_5} to multiple cells, with possibly common membranes (a.k.a. \textit{gap junctions}), is straightforward \cite{huynh2022convergence}.

\section{Weak formulation and discrete operators}\label{sec:approx}

The EMI problem can be weakly formulated in various ways, depending on the unknowns of interest. We refer to \cite{EMI_book} for a broad discussion on various formulations (including the so-called mixed ones). As it could be expected from the structure of \eqref{eq::EMI}-\eqref{eq::EMI_5}, all formulations and corresponding discretizations give rise to block operators, with different blocks corresponding to $\Omega_i,\Omega_e$, and possibly $\Gamma$.

We use the so-called \textit{single-dimensional} formulation and the corresponding discrete operators. In this setting, the weak form depends only on bulk quantities $u_i$ and $u_e$ since the current term $I_m$ is replaced by:
$$ I_m(\xx) = \tau^{-1}(u_i(\xx)-u_e(\xx)-f(\xx)),$$
according to equations~\eqref{eq::EMI_4}-\eqref{eq::EMI_5}. Let us remark that ``\textit{single}'' refers to the previous substitution, eliminating the variable defined in $\Gamma$; the overall EMI problem is still in multiple dimensions, in the sense that $d\geq1$.

After substituting the expression for $I_m$ in \eqref{eq::EMI_3}, assuming the solution $u_r$, for $r\in\{i,e\}$, to be sufficiently regular over $\Omega_r$, we multiply the PDEs in \eqref{eq::EMI}-\eqref{eq::EMI_2} by test functions $v_r(\xx)\in V_r(\Omega_r)$, with $V_r$ a sufficiently regular Hilbert space with elements satisfying the boundary conditions in \eqref{bc1}-\eqref{bc2}; in practice $H^1(\Omega_i)$ and $H^1(\Omega_e)$, would be a standard choice. After integrating over $\Omega_r$ and applying integration by parts, using the normal flux definition \eqref{eq::EMI_3}, the weak EMI problem reads: find $u_i\in V_i(\Omega_i)$ and $u_e\in V_e(\Omega_e)$ such that

\begin{align}
     \tau\int_{\Omega_e} \sigma_e\nabla u_e \cdot\nabla v_e\,\mathrm{d}\xx + \int_{\Gamma}u_ev_e\,\mathrm{d}s - \int_{\Gamma}u_iv_e\,\mathrm{d}s & = -\int_{\Gamma}fv_e\,\mathrm{d}s, \label{eq::EMI_weak1} \\
     \tau\int_{\Omega_i} \sigma_i\nabla u_i \cdot\nabla v_i\,\mathrm{d}\xx + \int_{\Gamma}u_iv_i\,\mathrm{d}s - \int_{\Gamma}u_ev_i\,\mathrm{d}s & = \int_{\Gamma}fv_i\,\mathrm{d}s, \label{eq::EMI_weak2}
\end{align}
for all test functions $v_e\in V_e(\Omega_e)$ and $v_i\in V_i(\Omega_i)$. We refer to \cite[Section 6.2.1]{EMI_book} for boundedness and coercivity results for this formulation. 

\begin{figure}   
    \centering        
    \begin{tikzpicture}

      \draw  plot[very thick, smooth, tension=.9] coordinates {(-6,0.5) (-6,3.5) (-4,2.5) (-3,3.5) (-2,1) (-5,0) (-6,0.5)};

  \draw[rounded corners=15pt]
  (-4,0.5) rectangle ++(1,1.5);

    \draw[rounded corners=15pt]
  (-5.5,1.3) rectangle ++(1.2,1.2);

    \node at (-2.8,2.4) {$\Omega_e$};
    \node at (-4.9,1.9) {$\Omega_i^1$};
    \node at (-3.5,1.2) {$\Omega_i^2$};
    \node at (-4.3,3.2) {$\partial\Omega$};
    \node at (-4.5,0.9) {$\Gamma$};
    
    \draw[] (-4,2.5) -- (-4.2,3);
    \draw[] (-5,1.3) -- (-4.7,1);
    \draw[] (-4,0.9) -- (-4.3,0.9);

    \draw[draw=black,very thick] (0,0) rectangle (4,4);
    \draw[draw=black,very thick] (1,1) rectangle (3,3);

    \node at (1.8,0.4) {$V_{e,h}(\Omega_e)$};
    \node at (1.8,2.4) {$V_{i,h}(\Omega_i)$};
    \node at (0.6,1.5) {$\Gamma$};
    \node at (-0.5,2.5) {$\partial\Omega$};

    \draw[] (0.8,1.5) -- (1,1.5);
    \draw[] (-0.2,2.5) -- (0,2.5);

    \draw[gray, very thin] (0,0) -- (4,4);
    \draw[gray, very thin] (1,0) -- (4,3);
    \draw[gray, very thin] (2,0) -- (4,2);
    \draw[gray, very thin] (3,0) -- (4,1);
    \draw[gray, very thin] (0,1) -- (3,4);
    \draw[gray, very thin] (0,2) -- (2,4);
    \draw[gray, very thin] (0,3) -- (1,4);

    \draw[gray, very thin] (1,0) -- (1,4);
    \draw[gray, very thin] (2,0) -- (2,4);
    \draw[gray, very thin] (3,0) -- (3,4);

    \draw[gray, very thin] (0,1) -- (4,1);
    \draw[gray, very thin] (0,2) -- (4,2);
    \draw[gray, very thin] (0,3) -- (4,3);

    \filldraw [red] (0,0) circle (2pt);
    \filldraw [red] (0,1) circle (2pt);
    \filldraw [red] (0,2) circle (2pt);
    \filldraw [red] (0,3) circle (2pt);
    \filldraw [red] (0,4) circle (2pt);
    \filldraw [red] (4,0) circle (2pt);
    \filldraw [red] (4,1) circle (2pt);
    \filldraw [red] (4,2) circle (2pt);
    \filldraw [red] (4,3) circle (2pt);
    \filldraw [red] (4,4) circle (2pt);
    \filldraw [red] (1,0) circle (2pt);
    \filldraw [red] (2,0) circle (2pt);
    \filldraw [red] (3,0) circle (2pt);
    \filldraw [red] (1,4) circle (2pt);
    \filldraw [red] (2,4) circle (2pt);
    \filldraw [red] (3,4) circle (2pt);

    \filldraw [cyan] (1.05,1.05) circle (2pt);
    \filldraw [cyan] (1.05,2) circle (2pt);
    \filldraw [cyan] (1.05,2.95) circle (2pt);
    \filldraw [cyan] (2.95,1.05) circle (2pt);
    \filldraw [cyan] (2.95,2) circle (2pt);
    \filldraw [cyan] (2.95,2.95) circle (2pt);
    \filldraw [cyan] (2,1.05) circle (2pt);
    \filldraw [cyan] (2,2.95) circle (2pt);

    \filldraw [orange] (0.95,0.95) circle (2pt);
    \filldraw [orange] (0.95,2) circle (2pt);
    \filldraw [orange] (0.95,3.05) circle (2pt);
    \filldraw [orange] (3.05,0.95) circle (2pt);
    \filldraw [orange] (3.05,2) circle (2pt);
    \filldraw [orange] (3.05,3.05) circle (2pt);
    \filldraw [orange] (2,0.95) circle (2pt);
    \filldraw [orange] (2,3.05) circle (2pt);

    \filldraw [blue] (2,2) circle (2pt);

    \filldraw [red] (4.7,3.1) circle (2pt);
    \filldraw [orange] (4.7,2.4) circle (2pt);
    \filldraw [blue] (4.7,1.7) circle (2pt);
    \filldraw [cyan] (4.7,1) circle (2pt);

    \node at (5.2,3.1) {$e$,in};
    \node at (5.2,2.4) {$e,\Gamma$};
    \node at (5.2,1.7) {$i$,in};
    \node at (5.2,1)   {$i,\Gamma$};

    \end{tikzpicture}
    \caption{Left: example of EMI geometry with $\Omega_i=\Omega_i^1\cup\Omega_i^2$. Right: example of the discretization setting for $d=2$ on a regular grid. We notice that interior and exterior nodes overlap in $\Gamma$. For $p = 1$, in this case, we have $N_i = 9, N_e = 24$ with $N_{i,\text{in}}= 1, N_{e,\text{in}} = 16$, and $N_\Gamma=8$. Corresponding labels are reported.}
    \label{fig:EMI}
\end{figure}
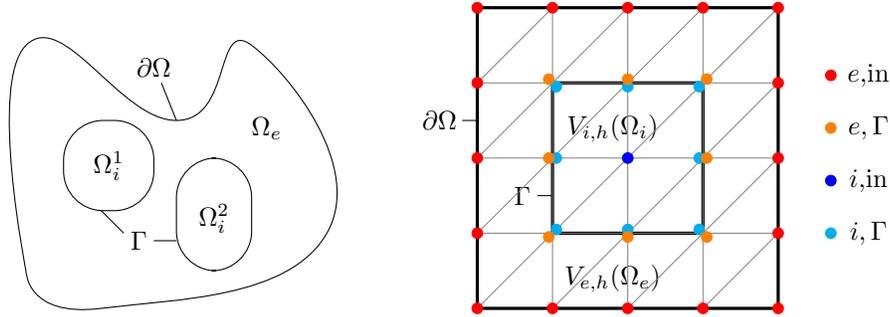

For each subdomain $\Omega_r$, with $r\in\{i,e\}$, we construct a conforming tassellation $\mathcal{T}_r$. We then introduce a yet unspecified discretization via finite element basis functions (e.g. Lagrangian elements of order $p\in\mathbb{N}$ on a regular grid) for $V_e$ and $V_i$:
$$V_{e,h}=\mathrm{span}\left(\{\phi^e_j\}_{j=1}^{N_e}\right),\quad V_{i,h}=\mathrm{span}\left(\{\phi^i_j\}_{j=1}^{N_i}\right),$$
with $N_e,N_i\in\mathbb{N}$ denoting the number of degrees of freedom in the corresponding subdomains. We can further decompose $N_r=N_{r,\text{in}}+N_{r,\Gamma}$, i.e. further differentiating between internal and membrane degrees of freedom, with $N_{r,\Gamma}$ basis functions $\phi^r_j$ having support intersecting $\Gamma$, 
cf. Figure~\ref{fig:EMI}. In the numerical experiments, we will consider matching $\mathcal{T}_r$ on the interface $\Gamma$ and the same $p$ for $V_{e,h}$ and $V_{i,h}$; nevertheless Section~\ref{sec:approx} and Section \ref{sec:spectral} are developed in a general setting, so that the theory is ready also for potential extensions.

%
From \eqref{eq::EMI_weak1}-\eqref{eq::EMI_weak2} we define the following discrete operators: intra- and extra- Laplacians
\begin{align}
    & A_e=\left[\int_{\Omega_e}\sigma_e\nabla\phi^e_j(\xx)\cdot\nabla\phi^e_k(\xx)\,\mathrm{d}\xx\right]_{j,k=1}^{N_e}\in\mathbb{R}^{N_e\times N_e}, \\
    & A_i=\left[\int_{\Omega_i}\sigma_i\nabla\phi^i_j(\xx)\cdot\nabla\phi^i_k(\xx)\,\mathrm{d}\xx\right]_{j,k=1}^{N_i}\in\mathbb{R}^{N_i\times N_i},
\end{align}
membrane mass matrices:
\begin{align}
    & M_{e}=\left[\int_{\Gamma}\phi^e_j(\xx)\phi^e_k(\xx)\,\mathrm{d}s\right]_{j,k=1}^{N_e}\in\mathbb{R}^{N_e\times N_e}, \\
    & M_{i}=\left[\int_{\Gamma}\phi^i_j(\xx)\phi^i_k(\xx)\,\mathrm{d}s\right]_{j,k=1}^{N_i}\in\mathbb{R}^{N_i\times N_i},
\end{align}
and the coupling matrix
\begin{align}
    & T_{ei}=\left[\int_{\Gamma}\phi^i_j(\xx)\phi^e_k(\xx)\,\mathrm{d}s\right]_{(j,k)=(1,1)}^{(N_e,N_i)}\in\mathbb{R}^{N_e\times N_i}.
\end{align}
Finally, we write the linear system of size $n = N_e+N_i$ corresponding to \eqref{eq::EMI_weak1}-\eqref{eq::EMI_weak2} as
\begin{equation}\label{eq:system}
	\begin{bmatrix}
		\tau A_e + M_e & T_{ei} \\
		T_{ei}^T & \tau A_i + M_i \\
	\end{bmatrix}
	\begin{bmatrix}
		\uu_e \\
		\uu_i \\
	\end{bmatrix} =
	\begin{bmatrix}
		\ff_e \\
		\ff_i \\
	\end{bmatrix} \quad \Longleftrightarrow \quad A_n \uu = \ff,
\end{equation}
with $\uu_r,\ff_r\in\mathbb{R}^{N_r}$ the unknowns and the right hand side corresponding to $\Omega_r$ and $r\in\{i,e\}$. We remark that the operator $A_n$ is symmetric and positive definite\footnote{Dirichlet boundary conditions can be imposed to the right-hand side to enforce symmetry.}. Making the degrees of freedom in $\Gamma$ denoted explicitly (with reference to Figure~\ref{fig:EMI}), i.e. $\uu_{i,\Gamma},\uu_{e,\Gamma}\in\mathbb{R}^{N_\Gamma}$, as well as the interior ones, i.e.  $\uu_{i,\mathrm{in}}\in\mathbb{R}^{N_i - N_\Gamma}$ and $\uu_{e,\mathrm{in}}\in\mathbb{R}^{N_e - N_\Gamma}$, we can rewrite more extensively system \eqref{eq:system} as
\begin{equation}
	\begin{bmatrix}
		\tau A_{e}^{\mathrm{in}} & \tau A_{e}^{\mathrm{in},\Gamma} & 0 & 0  \\
		\tau A_{e}^{\Gamma,\mathrm{in}} & \tau A_{e}^{\Gamma} + M_e^\Gamma & 0 & T_{ei}^\Gamma  \\
        0 & 0 & \tau A_{i}^{\mathrm{in}} & \tau A_{i}^{\mathrm{in},\Gamma}  \\
        0 & (T_{ei}^\Gamma)^T & \tau A_{i}^{\Gamma,\mathrm{in}} & \tau A_{i}^{\Gamma} + M_i^\Gamma  \\
		
	\end{bmatrix}
	\begin{bmatrix}
		\uu_{e,\mathrm{in}} \\
		\uu_{e,\Gamma} \\
		\uu_{i,\mathrm{in}} \\
		\uu_{i,\Gamma} \\
	\end{bmatrix} =
	\begin{bmatrix}
		\mathbf{0} \\
		\ff_{\Gamma,e} \\
		\mathbf{0}\\
		\ff_{\Gamma,i} \\
	\end{bmatrix},
\end{equation}
with
\begin{equation}\label{eq:system_matrix}
A_n=
\begin{bmatrix}
		\tau A_{e}^{\mathrm{in}} & \tau A_{e}^{\mathrm{in},\Gamma} & 0 & 0  \\
		\tau A_{e}^{\Gamma,\mathrm{in}} & \tau A_{e}^{\Gamma} + M_e^\Gamma & 0 & T_{ei}^\Gamma  \\
        0 & 0 & \tau A_{i}^{\mathrm{in}} & \tau A_{i}^{\mathrm{in},\Gamma}  \\
        0 & (T_{ei}^\Gamma)^T & \tau A_{i}^{\Gamma,\mathrm{in}} & \tau A_{i}^{\Gamma} + M_i^\Gamma  \\
			\end{bmatrix}.
\end{equation}
For forthcoming use, we define the bulk matrix
\[
\tilde A_n=\begin{bmatrix}
		\tau A_{e}^{\mathrm{in}} & 0 & 0 & 0  \\
		0 & 0 & 0 & 0  \\
        0 & 0 & \tau A_{i}^{\mathrm{in}} & 0  \\
        0 & 0 & 0 & 0  \\
			\end{bmatrix},
\]
where all membrane terms are zeroed.

\section{Spectral analysis}\label{sec:spectral}

In this section, we study the spectral distribution of the matrix sequence $\{A_{n}\}_n$ under various assumptions for determining the global behaviour of the eigenvalues of $A_n$ as the matrix size $n$ tends to infinity. The spectral distribution is given by a smooth function called the (spectral) symbol as it is customary in the Toeplitz and Generalized Locally Toeplitz (GLT) setting  \cite{GS-I,GS-II,block-glt-1D,block-glt-dD}.

First, we give the formal definition of Toeplitz structures, eigenvalue (spectral) and singular value distribution, the basic tools that we use from the GLT theory, and finally, we provide the specific analysis of our matrix sequences under a variety of assumptions.

\subsection{Toeplitz structures, spectral symbol, GLT tools}\label{ssec:GLTbackground}

We initially formalize the definition of block Toeplitz and circulant sequences associated with a matrix-valued Lebesgue integrable function. Then, we provide the notion of eigenvalue (spectral) and singular value distribution, and we introduce the basic tools taken from the GLT theory.


\begin{definition}\label{def_not1}{\rm [Toeplitz, block-Toeplitz, multilevel Toeplitz matrices]}
A finite-dimensional or infinite-dimensional Toeplitz matrix is a matrix that has constant
elements along each descending diagonal  from left to right, namely,
\begin{equation}\label{eq_toeplits}
T =
\begin{bmatrix}
a_0&a_{-1}&a_{-2}&\cdots\\
a_{1}&a_0&a_{-1}&\cdots\\
a_{2}&a_{1}&\ddots &\ddots\\
\ddots&\ddots&\ddots&\ddots
\end{bmatrix},
\, 
T_n =
\begin{bmatrix}
a_0&a_{-1}&a_{-2}&\cdots&a_{-n+1}\\
a_{1}&a_0&a_{-1}&\ddots&\\
a_{2}&a_{1}&\ddots&\ddots&\vdots\\
\vdots&&a_{1}&a_0&a_{-1}\\
a_{n-1}&\cdots&a_2&a_1&a_0
\end{bmatrix}
\end{equation}
As the first equation in \eqref{eq_toeplits} indicates, the matrix $T$ can be infinite-dimensional. A finite-dimensional Toeplitz matrix of dimension $n$ is denoted
as $T_n$.
We also consider sequences of Toeplitz matrices as a function of their
dimension, denoted as $\{T_n\}_n,\, n=1,2,\cdots,\infty$.

\smallskip
In general, the entries  $a_k$, $k\in \ZZ$, can be matrices (blocks) themselves, defining
$T_n$ as a block Toeplitz matrix. 
Thus, in the block case, 
$a_k, k\in \{-n+1, \cdots, n-1\}$ are blocks of size $s_1\times s_2$,
the subindex of $T_n$ is the number of blocks in the
Toeplitz matrix while $N_1\times N_2$ is the size of the 
matrix  with $N_1 = n\,s_1$, $N_2 = n\,s_2$.
To be specific, we use also the notation $X_{N_1,N_2}=T_n$. 
A
special case of block-Toeplitz matrices is the class of two- and multilevel
block Toeplitz matrices, where the blocks are Toeplitz (or multilevel
Toeplitz) matrices themselves. The standard Toeplitz matrices are sometimes
addressed as \textit{unilevel} Toeplitz.
\end{definition}

\begin{definition}\label{def_not3}{\rm [Toeplitz sequences (generating function of)]}
Denote by $f$ a $d$-variate complex-valued integrable function,
defined over the domain $Q^d=[-\pi,\pi]^d, d\ge 1$, with $d$-dimensional Lebesgue measure $\mu_d (Q^d)=(2\pi)^d$. Denote by $f_k$ the Fourier coefficients of $f$,
$$
f_k = \frac{1}{(2\pi)^d}\int\limits_{Q^d}f(\theta)e^{-i\,(k,\theta)}\,d\theta,\;
k=(k_1,\cdots,k_d)\in \ZZ^d, \; i^2=-1,
$$
where $\theta=(\theta_1,\cdots,\theta_d)$, $(k,{ \theta})=\sum_{j=1}^{d}k_j \theta_j$, $n=(n_1,\cdots,n_d)$, and 
$N(n)={n_1}\times\cdots\times n_d$. By following the multi-index notation in \cite{tyrt}[Section 6], with each
$f$ we can associate a sequence of Toeplitz matrices $\{T_n\}_n$, where
$$
T_n = \{ f_{k-\ell}\}_{k,\ell=\mathbf{e}^T}^{n} \in \mathbb{C}^{N(n)\times N(n)},
$$
$\mathbf{e}=[1,1,\cdots,1]\in \NN^d$.
For $d=1$
$$
T_n =
\begin{bmatrix}
f_0&f_{-1}&\cdots&f_{-n+1}\\
f_{1}&f_0&\ddots&\vdots\\
\vdots&\ddots&f_0&f_{-1}\\
f_{n-1}&\cdots&f_1&f_0
\end{bmatrix},
$$
or for $d = 2$, i.e. the two-level case, and for example $n=(2,3)$, we have
$$
T_n =
\begin{bmatrix}
F_0&F_{-1}\\
F_{1}&F_0
\end{bmatrix},
\quad
F_k =
\begin{bmatrix}
f_{k,0}&f_{k,-1}&f_{k,-2}\\
f_{k,1}&f_{k,0}&f_{k,-1}\\
f_{k,2}&f_{k,1}&f_{k,0}
\end{bmatrix},\ \ \ k=0,\pm 1.
 $$
The function $f$ is referred to as the \textit{generating function} (or the
\textit{symbol} of) $T_n$.  Using a more compact notation, we say that the
function $f$ is the generating function of the whole sequence $\{T_n\}_n$ and we
write $T_n = T_n(f)$. 

If $f$ is $d$-variate, $\mathbb{C}^{s_1\times s_2}$ matrix-valued, and integrable over  $Q^d$, $d,s_1,s_2\ge 1$, i.e. 
$f\in L^1(Q^2,s_1\times s_2)$, then we can define the
Fourier coefficients of $f$ in the same way (now $f_k$ is a matrix of size
$s_1\times s_2$) and consequently  $T_n = \{ f_{k-\ell}\}_{k,\ell=\mathbf{e}^T}^{n}
\in \mathbb{C}^{s_1N(n)\times s_2N(n)}$, then $T_n$ is a $d$-level block
Toeplitz matrix according to Definition \ref{def_not1}. 
If $s_1=s_2=s$ then we write $f\in L^1(Q^2,s)$.

As in the scalar case, the function $f$ is referred to as the \textit{generating function} of $T_n$. We say that the function $f$ is the generating
function of the whole sequence $\{T_n\}_n$, and we use the notation $T_n = T_n(f)$. 
\end{definition}

\begin{definition}\label{def-distribution}
	Let $f:D\to\mathbb{C}^{s\times s}$ be a measurable matrix-valued function with eigenvalues $\lambda_i(f)$ and singular values $\sigma_i(f)$, $i=1,\ldots,s$. Assume that $D\subset \mathbb{R}^d$ is Lebesgue measurable with positive and finite Lebesgue measure $\mu_d(D)$. Assume that $\{A_n\}_n$ is a sequence of matrices such that ${\rm dim}(A_n)=d_n\rightarrow\infty$, as $n\rightarrow\infty$ and with eigenvalues $\lambda_j(A_n)$ and singular values $\sigma_j(A_n)$, $j=1,\ldots,d_n$.
	\begin{itemize}
		\item We say that $\{A_n\}_{n}$ is {\em distributed as $f$ over $D$ in the sense of the eigenvalues,} and we write $\{A_n\}_{n}\sim_\lambda(f,D),$ if
		\begin{equation}\label{distribution:eig}
		\lim_{n\to\infty}\frac{1}{d_n}\sum_{j=1}^{d_n}F(\lambda_j(A_n))=
		\frac1{\mu_d(D)} \int_D \frac1{s}\sum_{i=1}^sF(\lambda_i(f(t))) \, \mathrm{d}t,
		\end{equation}
		for every continuous function $F$ with compact support. In this case, we say that $f$ is the \emph{spectral symbol} of $\{A_{n}\}_{n}$.
		
		\item We say that $\{A_n\}_{n}$ is {\em distributed as $f$ over $D$ in the sense of the singular values,} and we write $\{A_n\}_{n}\sim_\sigma(f,D)$, if
		\begin{equation}\label{distribution:sv}
		\lim_{n\to\infty}\frac{1}{d_n}\sum_{j=1}^{d_n}F(\sigma_j(A_n))=
		\frac1{\mu_d(D)} \int_D \frac1{s}\sum_{i=1}^sF(\sigma_i(f(t))) \, \mathrm{d}t,
		\end{equation}
		for every continuous function $F$ with compact support. In this case, we say that $f$ is the \emph{singular value symbol} of $\{A_{n}\}_{n}$.
		\item The notion $\{A_n\}_{n}\sim_\sigma(f,D)$ applies also in the rectangular case where $f$ is $\mathbb{C}^{s_1\times s_2}$ matrix-valued. In such a case the parameter $s$ in formula (\ref{distribution:sv}) has to be replaced by the minimum between $s_1$ and $s_2$: furthermore $A_n\in \mathbb{C}^{d_n^{(1)}\times d_n^{(2)}}$ with $d_n$ in formula (\ref{distribution:sv}) being the minimum between $d_n^{(1)}$ and $d_n^{(2)}$. Of course the notion of eigenvalue distribution does not apply in a rectangular setting.
			\end{itemize}
\end{definition}
Throughout the paper, when the domain can be easily inferred from the context, we replace the notation $\{A_n\}_n\sim_{\lambda,\sigma}(f,D)$ with $\{A_n\}_n\sim_{\lambda,\sigma} f$.

\begin{remark}\label{rem:approx}
	If $f$ is smooth enough, an informal interpretation of the limit relation
	\eqref{distribution:eig} (resp. \eqref{distribution:sv})
	is that when $n$ is sufficiently large, the eigenvalues (resp. singular values) of $A_{n}$ can be subdivided into $s$ different subsets of the same cardinality. Then $d_n/s$ eigenvalues (resp. singular values) of $A_{n}$ can
	be approximated by a sampling of $\lambda_1(f)$ (resp. $\sigma_1(f)$)
	on a uniform equispaced grid of the domain $D$, and so on until the
	last $d_n/s$ eigenvalues (resp. singular values), which can be approximated by an equispaced sampling
	of $\lambda_s(f)$ (resp. $\sigma_s(f)$) in the domain $D$.
\end{remark}
\begin{remark}
   We say that $A_n$ is \textit{zero-distributed} if $A_n\sim_\lambda 0$. Of course, if the eigenvalues of $A_n$ tend all to zero for $n\to\infty$, then this is sufficient to claim that $A_n\sim_\lambda 0$.
\end{remark}
For Toeplitz matrix sequences, the following theorem due to Tilli holds, which generalizes previous research along the last 100 years by Szeg\H{o}, Widom, Avram, Parter, Tyrtyshnikov, and Zamarashkin (see \cite{GS-II,block-glt-dD} and references therein).

\begin{theorem}{\rm \cite{TilliNota}}\label{szego-herm}
	Let $f\in L^1(Q^d,s_1\times s_2)$, then $\{T_{n}(f)\}_{{n}}\sim_\sigma(f,Q^d).$ If $s_1=s_2=s$ and if $f$ is a Hermitian matrix-valued function, then $\{T_{n}(f)\}_{{n}}\sim_\lambda(f,Q^d)$.
\end{theorem}

The following theorem is useful for computing the spectral distribution of a sequence of Hermitian matrices. For the related proof, see \cite[Theorem 4.3]{curl-div} and \cite[Theorem 8]{NS1-elonged}. Here, the conjugate transpose of the matrix $X$ is denoted by $X^*$.

\begin{theorem}{\rm \cite[Theorem 4.3]{curl-div}}\label{th:extradimensional}
	Let $\{A_n\}_n$ be a sequence of matrices, with $A_n$ Hermitian of size $d_n$, and let $\{P_n\}_n$ be a sequence such that $P_n\in\mathbb C^{d_n\times\delta_n}$, $P_n^*P_n=I_{\delta_n}$, $\delta_n\le d_n$ and $\delta_n/d_n\to1$ as $n\to\infty$. Then $\{A_n\}_n\sim_{\lambda}f$ if and only if $\{P_n^*A_nP_n\}_n\sim_{\lambda}f$.
\end{theorem}

With the notations of the result above, the matrix sequence $\{P_n^*A_nP_n\}_n$ is called a compression of $\{A_n\}_n$ and the single matrix $P_n^*A_nP_n$ is called a compression of $A_n$.

In what follows we take into account a crucial fact that often is neglected: the generating function of a Toeplitz matrix sequence and even more the spectral symbol of a given matrix sequence is not unique, except for the trivial case of either a constant generating function or a constant spectral symbol. In fact, here we report and generalize Remark 1.3 at p. 76 in \cite{Schur-CMAME} and the discussion below Theorem 3 at p. 8 in \cite{Symbols-SE}.

\begin{remark}\label{rem:non-uniq}{\rm [Multiple Toeplitz generating functions and multiple spectral symbols]}
Let $n$ be a multiple of $k$ and consider the Toeplitz matrix $X_n = T_n(f)$. We can view $X_n$ as a block-Toeplitz with blocks of size $k \times k$ such that $X_n = T_n(f) = T_{\frac{n}{k}}(f^{[k]})$. 
In our specific context, we have
\[
f^{[k]}(\theta)= T_k(f)-{\mathbf e}_1{\mathbf e}_k^T e^{i\ \theta} -{\mathbf e}_k {\mathbf e}_1^T e^{-i\ \theta}
\]
with ${\mathbf e}_j$, $j=1,\ldots,k$, being the canonical basis of $\mathbb{C}^k$.

As an example, let $n$ be even and consider the function $f(\theta)=2-2\cos(\theta)$. According to Definition \ref{def_not3}, $X_n=T_n(f)$ but we can also view the matrix as $X_n = T_{\frac{n}{2}}(f^{[2]})$, namely
\[
X_n=T_n(f)=\begin{bmatrix}
     2          & -1     &        & \mathbf{0}\\
    -1          & \ddots & \ddots \\
                & \ddots & \ddots & -1\\
    \mathbf{0}  &        & -1     & 2
     \end{bmatrix} = 
\begin{bmatrix}
A_0        & A_{-1} &        & \mathbf{0} \\
A_{1}      & \ddots & \ddots\\
           & \ddots & \ddots & A_{-1}\\
\mathbf{0} &        &        A_1 &A_0
\end{bmatrix} = T_{n\over 2}(f^{[2]}),
\]
where
\[
A_0=\begin{bmatrix} 2 & -1 \\ -1 & 2 \end{bmatrix},\ \ \ 
A_1=A_{-1}^T= \begin{bmatrix} 0 & -1 \\ 0 & 0 \end{bmatrix}.
\]
Thus, $f^{[2]} = \left(A_0+A_1 e^{i\theta} +A_1^T e^{-i\theta}\right)$.  

It is clear that analogous multiple representations of Toeplitz matrices hold
for every function $f$. As a consequence, taking into consideration Theorem \ref{szego-herm}, we have simultaneously $\{T_{n}(f)\}_{{n}}\sim_\lambda(f,Q^d)$ and 
$\{T_{n}(f)\}_{{n}}\sim_\lambda(f^{[k]} ,Q^d)$ for any fixed $k$ independent of $n$. It should be also observed that a situation in which $k$ does not divide $n$ is not an issue thanks to the compression argument in Theorem \ref{th:extradimensional}.

More generally, we can safely claim that the spectral symbol in the Weyl sense of Definition \ref{def-distribution} is far from unique and in fact any rearrangement is still a symbol (see \cite{Symbols-SE,barbarino2022constructive}). A simple case is given by standard Toeplitz sequences $\{T_n(f)\}_n$, with $f$ real-valued and  even that is $f(\theta)=f(-\theta)$ almost everywhere, $\theta\in Q$. In that case 
\begin{equation}\label{distribution:eig-T}
 {1\over 2\pi}\int_{-\pi}^{\pi} F(f(\theta))\,\mathrm{d}\theta=
{1\over \pi}\int_{0}^{\pi} F(f(\theta))\,\mathrm{d}\theta,
\end{equation}
due to the even character of $f$,  and hence it is also true that $\{T_n(f)\}_n\sim_{\lambda} (f,Q_+)$, $Q_+=(0,\pi)$. A general analysis of these concepts via rearrangement theory can be found in \cite{barbarino2022constructive} and references therein. 
\end{remark}

\subsection{Symbol analysis}\label{ssec:symbol_analysis}

In this section, we state and prove three results which hold under reasonable assumptions. The first has the maximal generality, while the second and the third can be viewed as special cases of the first. Let us remind that $N_i,N_e$ and $N_\Gamma$ denote the number of intra-, extra- and membrane degrees of freedom.  

\begin{theorem}\label{th:general}
Assume that
\begin{eqnarray*}
N_\Gamma  =  o(\min\{N_i,N_e\}) & &  \mathrm{for} \quad N_\Gamma,N_i,N_e\to\infty, \\
\lim_{N_i,N_e\rightarrow\infty}  \frac{N_i}{N_e+N_i} = r \in (0,1).
\end{eqnarray*}
Assume that
\begin{align}\label{distrAe}
    & \{\tau A_{e}^{\mathrm{in}}\}_n\sim_\lambda (f^e,D^e), \\ \label{distrAi}
    & \{\tau A_{i}^{\mathrm{in}}\}_n\sim_\lambda (f^i,D^i),
\end{align}
$D^e\subset \mathbb{R}^{k_e}$, $D^i\subset \mathbb{R}^{k_i}$, given $f^e$ and $f^i$ extra- and intra- spectral symbols.
It follows that
\[
\{\tilde A_n\}_n,\ \ \{A_n\}_n\sim_\lambda (g,[0,1]\times D^i\times D^e) 
\]
where $g(x,t_i,t_e)=f^i(t_i) \psi_{[0,r]}(x) + f^e(t_e) \psi_{(r,1]}(x)$, $x\in [0,1]$, $t_i\in D^i$, $t_e\in D^e$, 
and
\[
	\begin{bmatrix}
		0 & \tau A_{e}^{\mathrm{in},\Gamma} & 0 & 0  \\
		\tau A_{e}^{\Gamma,\mathrm{in}} & \tau A_{e}^{\Gamma} + M_e^\Gamma & 0 & T_{ei}^\Gamma  \\
        0 & 0 & 0 & \tau A_{i}^{\mathrm{in},\Gamma}  \\
        0 & (T_{ei}^\Gamma)^T & \tau A_{i}^{\Gamma,\mathrm{in}} & \tau A_{i}^{\Gamma} + M_i^\Gamma  \\
			\end{bmatrix}  \sim_\lambda 0,
\]
$\psi_Z$ denoting the characteristic function of the set $Z$.
\end{theorem}
\ \\
{\bf Proof}

First, we observe that because of the Galerkin approach and the structure of the equations, all the involved matrices are real and symmetric. Hence the symbols $f^e$, $f^i$ are necessarily Hermitian matrix-valued. Without loss of generality we assume that both $f^e$, $f^i$ take values into 
$\mathbb{C}^{s\times s}$ (the case where $f^e$, $f^i$ have different matrix sizes is treated in Remark \ref{rem-different s}, for the sake of notational simplicity).
Because of equations (\ref{distrAe}) and (\ref{distrAi}), setting $N_e'=N_e-N_\Gamma$, $N_i'=N_i-N_\Gamma$, and taking into account Definition \ref{def-distribution}, we have
\begin{eqnarray} \label{distr-e}
\lim_{N_e\to\infty}\frac{1}{N_e'}\sum_{j=1}^{N_e'}F(\lambda_j(\tau A_{e}^{\mathrm{in}})) & = & 
		\frac1{\mu_{k_e}(D^e)} \int_{D^e} \frac1{s}\sum_{m=1}^sF(\lambda_m(f^e(t_e)))\,\mathrm{d}t_e,  \\ \label{distr-i}
\lim_{N_i\to\infty}\frac{1}{N_i'}\sum_{j=1}^{N_i'}F(\lambda_j(\tau A_{i}^{\mathrm{in}})) & = & 
		\frac1{\mu_{k_i}(D^i)} \int_{D^i} \frac1{s}\sum_{m=1}^sF(\lambda_m(f^i(t_i))) \,\mathrm{d}t_i,		
\end{eqnarray}
for any continuous function $F$ with bounded support. 
Now we build the $2\times 2$ block matrix 
\[
X[N_e,N_i,N_\Gamma]={\rm diag}(\tau A_{e}^{\mathrm{in}},\tau A_{i}^{\mathrm{in}})
\]
and we consider again a generic continuous function $F$ with bounded support. By the block diagonal structure of $X[N_e,N_i,N_\Gamma]$, we infer 

\begin{eqnarray*}
\Delta(N_e,N_i,N_\Gamma,F) & = &  \frac{1}{N_e'+N_i'}\sum_{j=1}^{N_e'+N_i'}F(\lambda_j(X[N_e,N_i,N_\Gamma])) \\
 & = & \sum_{j=1}^{N_e'}\frac{F(\lambda_j(\tau A_{e}^{\mathrm{in}}))}{N_e'+N_i'} + \sum_{j=1}^{N_i'}\frac{F(\lambda_j(\tau A_{i}^{\mathrm{in}}))}{N_e'+N_i'} \\
 & = & \frac{N_e'}{N_e'+N_i'}\sum_{j=1}^{N_e'}\frac{F(\lambda_j(\tau A_{e}^{\mathrm{in}}))}{N_e'} + \frac{N_i'}{N_e'+N_i'}\sum_{j=1}^{N_i'}\frac{F(\lambda_j(\tau A_{i}^{\mathrm{in}}))}{N_i'}.
\end{eqnarray*}

As a consequence, taking the limit as $N_i,N_e\to \infty$, using the assumption $N_\Gamma  =  o(\min\{N_i,N_e\})$,
\[
\lim_{N_i,N_e\rightarrow\infty}  \frac{N_i}{N_e+N_i} = r \in (0,1),
\]
we obtain that the limit of $\Delta(N_e,N_i,N_\Gamma,F)$ exists and it is equal to
\begin{equation}\label{eq: ci-siamo}
\frac{r}{\mu_{k_e}(D^e)} \int_{D^e} \frac1{s}\sum_{m=1}^sF(\lambda_m(f^e(t_e))) \,\mathrm{d}t_e +
\frac{1-r}{\mu_{k_i}(D^i)} \int_{D^i} \frac1{s}\sum_{m=1}^sF(\lambda_m(f^i(t_i))) \,\mathrm{d}t_i.
\end{equation}
With regard to Definition \ref{def-distribution}, the quantity (\ref{eq: ci-siamo}) does not look like the right-hand side of (\ref{distribution:eig}), since we see two different terms. The difficulty can be overcome by enlarging the space with a fictitious domain $[0,1]$ and interpreting the sum in (\ref{eq: ci-siamo}) as the global integral involving a step function. In reality, we rewrite (\ref{eq: ci-siamo}) as
\[
\frac{1}{\mu_{q}(D)} 
\int_0^1 \,\mathrm{d}x \int_{D^e\times D^i} \frac1{s}\sum_{m=1}^sF(\lambda_m(f^e(t_e)) \psi_{[0,r]}(x) + \lambda_m(f^i(t_i)) \psi_{(r,1]}(x))\,\mathrm{d}t_e\mathrm{d}t_i,
\]
with $D=[0,1]\times D^i\times D^e$, $q=k_e+k_i+1$, $\mu_{q}(D)=\mu_{k_e}(D^e)\mu_{k_i}(D^i)$. Hence the proof of the relation
\[
\{X[N_e,N_i,N_\Gamma]\}_n\sim_\lambda (g,[0,1]\times D^i\times D^e) 
\]
is complete.
Now we observe that $X[N_e,N_i,N_\Gamma]$ is a compression of $\tilde A_n$ (cf. definition after Theorem~\ref{th:extradimensional}): therefore by exploiting again the hypothesis $N_\Gamma  =  o(\min\{N_i,N_e\})$ which implies  $N_\Gamma  =  o(n)$, by Theorem \ref{th:extradimensional}, we deduce  
\[
\{\tilde A_n\}_n \sim_\lambda (g,[0,1]\times D^i\times D^e). 
\]
Finally in Exercise 5.3 in the book \cite{GS-I} it is proven the following: if $\{V_n\}_n\sim_{\lambda}f$ (not necessarily Toeplitz or GLT), $\{A_n=V_n+Y_n\}_n$  with $\{Y_n\}_n$  zero-distributed and $A_n,Y_n$ Hermitian for every dimension, then $\{A_n\}_n\sim_{\lambda}f$.
This setting is exactly our setting since by direct inspection the rank of 
\[
R_n=	\begin{bmatrix}
		0 & \tau A_{e}^{\mathrm{in},\Gamma} & 0 & 0  \\
		\tau A_{e}^{\Gamma,\mathrm{in}} & \tau A_{e}^{\Gamma} + M_e^\Gamma & 0 & T_{ei}^\Gamma  \\
        0 & 0 & 0 & \tau A_{i}^{\mathrm{in},\Gamma}  \\
        0 & (T_{ei}^\Gamma)^T & \tau A_{i}^{\Gamma,\mathrm{in}} & \tau A_{i}^{\Gamma} + M_i^\Gamma  \\
			\end{bmatrix}  
\]
is bounded by $4N_\Gamma$. Hence the number of nonzero eigenvalues of the latter matrix cannot exceed $4N_\Gamma=o(n)$, $n=N_i+N_e$, and therefore the related matrix sequence is zero-distributed in the eigenvalue sense i.e. $\{R_n\}_n\sim_\lambda 0$. Since $A_n=\tilde A_n + R_n$, the application of the statement in \cite[Exercise 5.3]{GS-I} implies directly $\{A_n\}_n\sim_\lambda (g,[0,1]\times D^i\times D^e)$ and the proof is concluded. 
\ \hfill $\bullet$

\begin{remark}\label{rem-different s}
The case where $f^e$ and $f^i$ have different matrix sizes would complicate the derivations in the proof of the above theorem. 
However, in this case, the argument is simple and relies on the non-uniqueness of the spectral symbol: for a discussion on the matter refer to \cite{Symbols-SE,barbarino2022constructive} and Remark \ref{rem:non-uniq}.  
\end{remark}

The following two corollaries simplify the statement of Theorem \ref{th:general}, under special assumptions which are satisfied for a few basic discretization schemes and when dealing with elementary domains.

\begin{corollary}\label{th:specific}
Assume that
\begin{eqnarray*}
N_\Gamma = o(\min\{N_i,N_e\}) &  &  \mathrm{for} \quad N_\Gamma,N_i,N_e\to\infty, \\
\lim_{N_i,N_e\rightarrow\infty}  \frac{N_i}{N_e+N_i} = 1/2. 
\end{eqnarray*}
Assume that
\begin{align}
    & \{\tau A_{e}^{\mathrm{in}}\}_n\sim_\lambda (f^e,D), \\
    & \{\tau A_{i}^{\mathrm{in}}\}_n\sim_\lambda (f^i,D),
\end{align}
$D\subset \mathbb{R}^{k}$.
It follows that
\[
\{\tilde A_n\}_n,\ \ \{A_n\}_n\sim_\lambda \left(\begin{bmatrix}
		f^e & 0 \\
		0 & f^i \\
	\end{bmatrix},D\right)
\]
and
\[
	\begin{bmatrix}
		0 & \tau A_{e}^{\mathrm{in},\Gamma} & 0 & 0  \\
		\tau A_{e}^{\Gamma,\mathrm{in}} & \tau A_{e}^{\Gamma} + M_e^\Gamma & 0 & T_{ei}^\Gamma  \\
        0 & 0 & 0 & \tau A_{i}^{\mathrm{in},\Gamma}  \\
        0 & (T_{ei}^\Gamma)^T & \tau A_{i}^{\Gamma,\mathrm{in}} & \tau A_{i}^{\Gamma} + M_i^\Gamma  \\
			\end{bmatrix} \sim_\lambda 0.
\]
\end{corollary}
{\bf Proof}

Notice that the writing
\[
\{\tilde A_n\}_n,\ \ \{A_n\}_n\sim_\lambda \begin{bmatrix}
		f^e & 0 \\
		0 & f^i \\
	\end{bmatrix}
\]
and 
\[
\{\tilde A_n\}_n,\ \ \{A_n\}_n\sim_\lambda f^i(\theta) \psi_{[0,c]}(x) + f^e(\theta) \psi_{(c,1]}(x).
\]
are equivalent for $c=1/2$.

\ \hfill $\bullet$

\begin{corollary}\label{th:specific-bis}
Assume that
\begin{eqnarray*}
N_\Gamma = o(\min\{N_i,N_e\}) &  &  \mathrm{for} \quad N_\Gamma,N_i,N_e\to\infty, \\
\lim_{N_i,N_e\rightarrow\infty}  \frac{N_i}{N_e+N_i} = r\in [0,1]. 
\end{eqnarray*}
Assume that
\begin{align}
    & \{\tau A_{e}^{\mathrm{in}}\}_n\sim_\lambda (f,D), \\
    & \{\tau A_{i}^{\mathrm{in}}\}_n\sim_\lambda (f,D),
\end{align}
$D\subset \mathbb{R}^{k}$.
It follows that
\[
\{\tilde A_n\}_n,\ \ \{A_n\}_n\sim_\lambda (f,D)
\]
and
\[
	\begin{bmatrix}
		0 & \tau A_{e}^{\mathrm{in},\Gamma} & 0 & 0  \\
		\tau A_{e}^{\Gamma,\mathrm{in}} & \tau A_{e}^{\Gamma} + M_e^\Gamma & 0 & T_{ei}^\Gamma  \\
        0 & 0 & 0 & \tau A_{i}^{\mathrm{in},\Gamma}  \\
        0 & (T_{ei}^\Gamma)^T & \tau A_{i}^{\Gamma,\mathrm{in}} & \tau A_{i}^{\Gamma} + M_i^\Gamma  \\
			\end{bmatrix} \sim_\lambda 0.
\]
\end{corollary}
{\bf Proof}
The proof of the relation $\{\tilde A_n\}_n,\ \ \{A_n\}_n\sim_\lambda (f,D)$ follows directly from the limit displayed in (\ref{eq: ci-siamo}),
after replacing $f^e,f^i$ with $f$ and necessarily $D^e,D^i$ with $D$. The rest is obtained verbatim as in Theorem \ref{th:general}. 

\ \hfill $\bullet$

\subsection{Analysis of stiffness matrix sequences}\label{ssec:stiffness_matrices}
The spectral analysis of stiffness matrix sequences is crucial to understand the global spectral behaviour of EMI matrices and, in turn, to justify our preconditioning approach.
According to \cite[Section 5]{QpFEM}, given a bi-dimensional square domain, the sequence of corresponding stiffness matrices $\{A_n^{\rm{square}}\}_n$ obtained from the linear $\mathbb{Q}_1$ Lagrangian finite element method\footnote{The label $\mathbb{Q}_p$ denotes the space of elements of order $p$ with square support while $\mathbb{P}_p$ is used for elements in a triangulated mesh. For a broader understanding and the analysis with the $\mathbb{P}_1$ elements, we refer to \cite[Section 7]{reducedGLT}.} (FEM) has the following spectral distribution:
\[
	\{A_n^{\rm{square}}\}_n\sim_{\lambda} f^{\rm{square}},
\]
with $f^{\rm{square}}(\theta_1,\theta_2) = h_{\rm{1D}}(\theta_1)  f_{\rm{1D}}(\theta_2)  + f_{\rm{1D}}(\theta_1) h_{\rm{1D}}(\theta_2)$, the functions $f_{\rm{1D}}$ and $h_{\rm{1D}}$ being, respectively, the symbols associated to the stiffness and mass matrix sequences in one dimension, that is, $f_{\rm{1D}}(\theta) = 2-2\cos(\theta)$, $h_{\rm{1D}}(\theta) = \frac{2}{3}+\frac{1}{3}\cos(\theta)$. Expanding the coefficients, we have
\begin{equation}\label{eq:f_square}
   	f^{\rm{square}}(\theta_1,\theta_2) = \frac{8}{3}-\frac{2}{3}\cos(\theta_1)-\frac{2}{3}\cos(\theta_2)
						-\frac{2}{3}\cos{(\theta_1+\theta_2)}-\frac{2}{3}\cos(\theta_1-\theta_2). 
\end{equation}
Let us now consider $\{T_n(f^{\rm{square}})\}_n$, the sequence of bi-level Toeplitz matrices generated by $f^{\rm{square}}$, which is a multilevel GLT sequence by definition (see \cite{GS-II}). We adopt the bi-index notation, that is, we consider the elements of $T_n(f^{\rm{square}})$ to be indexed with a bi-index $(i_1,i_2)$ with $i_1 = 1,\dots,n_1$ and $i_2 = 1,\dots,n_2$.

Let us focus on $\Omega_i$. In the $\mathbb{Q}_1$ case there is a one-to-one correspondence between the mesh grid points which are in the interior of the subdomain $\Omega_i^\mathrm{in}$ and the $N_i$ degrees of freedom, where $N_i$ is the dimension of $V_{i,h}$. See \cite[Section 3]{QpFEM} for an explanation.

Following the analysis performed in \cite{reducedGLT}, first, we set to zero all the rows $(i_1,i_2)$ and columns $(i_1,i_2)$ in $T_n(f^{\rm{square}})$ such that $\left(\frac{i_1}{n_1+1},\frac{i_2}{n_2+1}\right)\not\in\Omega_i^\mathrm{in}$, which in the $\mathbb{Q}_1$ case means that $\left(\frac{i_1}{n_1+1},\frac{i_2}{n_2+1}\right)$ is not a mesh grid point tied to a degree of freedom associated with a trial function in $V_{i,h}$. We call the matrix resulting from this operation $T_i$. Then, we define the restriction maps $\Pi_{\Omega_i}$ and $\Pi_{\Omega_i}^T$ such that $\Pi_{\Omega_i} T_i \Pi_{\Omega_i}^T$ is the matrix obtained from $T_i$ deleting all rows $(i_1,i_2)$ and columns $(i_1,i_2)$ such that $\left(\frac{i_1}{n_1+1},\frac{i_2}{n_2+1}\right)\not\in\Omega_i^\mathrm{in}$. With a similar procedure, we construct also the matrices $\Pi_{\Omega_e} T_e \Pi_{\Omega_e}^T$. Then, we can apply \cite[Lemma 4.4]{reducedGLT} and conclude that 
\begin{align*}
    & \{\tau\Pi_{\Omega_e} T_e \Pi_{\Omega_e}^T\}_n\sim_\lambda (f^e,D^e), \\ 
    & \{\tau\Pi_{\Omega_i} T_i \Pi_{\Omega_i}^T\}_n\sim_\lambda (f^i,D^i),
\end{align*}
with $f^i(\xx,\boldsymbol{\theta}) = \tau f^{\rm{square}}(\boldsymbol{\theta})$, $(\xx,\boldsymbol{\theta}) \in D^i = \Omega_i \times [0,\pi]^2$, and $f^e(\xx,\boldsymbol{\theta}) = \tau f^{\rm{square}}(\boldsymbol{\theta})$, $(\xx,\boldsymbol{\theta}) \in D^e = \Omega_e \times [0,\pi]^2$. Apart from a low-rank correction due to boundary conditions, the matrices $\Pi_{\Omega_i} T_i \Pi_{\Omega_i}^T$ coincide with $A_{i}^{\mathrm{in}}$ and the same holds for the external domain. Since low-rank corrections are zero-distributed and hence in the Hermitian case have no impact on spectral distributions (see Exercise 5.3 in book \cite{GS-I}), we can safely write
\begin{align}\label{f^e-specific}
    & \{\tau A_{e}^{\mathrm{in}}\}_n\sim_\lambda (f^e,D^e), \\ \label{f^i-specific}
    & \{\tau A_{i}^{\mathrm{in}}\}_n\sim_\lambda (f^i,D^i).
\end{align}
Finally, by exploiting again the useful non-uniqueness of the spectral symbol as in Remark \ref{rem:non-uniq}, we also have 
\begin{align}\label{f^e-specific-bis}
    & \{\tau A_{e}^{\mathrm{in}}\}_n\sim_\lambda (\tau f^{\rm{square}},[0,\pi]^2), \\ \label{f^i-specific-bis}
    & \{\tau A_{i}^{\mathrm{in}}\}_n\sim_\lambda (\tau f^{\rm{square}},[0,\pi]^2),
\end{align}
since $f^{\rm{square}}$ over $[0,\pi]^2$ is a rearrangement of both $f^e$ over $D^e$ and $f^i$ over $D^i$.

Given relationships (\ref{f^e-specific}) and (\ref{f^i-specific}), Theorem \ref{th:general} applies since the other assumptions are verified when using any reasonable discretization. For instance, the parameter $r$ in Theorem \ref{th:general} is the ratio 
\[
\frac{\mu_d(\Omega_i)}{\mu_d(\Omega_i)+\mu_d(\Omega_e)}, \ \ \ d=1,2,3. 
\]
Furthermore, if $\mu_d(\Omega_i)=\mu_d(\Omega_e)$, $r=1/2$ and Corollary \ref{th:specific} applies. We observe that in our specific setting, due to (\ref{f^e-specific-bis}) and (\ref{f^i-specific-bis}), also the assumptions of Corollary \ref{th:specific-bis} hold true, so that the assumption on the ratio $r$ is no longer relevant.

\subsection{Monolithic solution strategy}\label{ssec:preconditioner}
Since, by construction, $A_n$ is symmetric positive definite (cf. Section~\ref{sec:approx}), we employ the conjugate gradient method (CG) to solve system \eqref{eq:system} iteratively. We first describe the dependency on $\tau$ and subsequently discuss preconditioning. 

\subsubsection*{Conjugate gradient convergence}
 To better understand the convergence properties of the CG method in the EMI context, in particular robustness w.r.t. $\tau$, we recall the following result, proved in \cite{axelsson1986rate}.

\begin{lemma}\label{lemma}
Assume $A_n$ to be a symmetric positive definite matrix of size $n$. Moreover, suppose there exists $a,b>0$, both independent of the matrix size $n$, and an integer $q<n$ such that 
\begin{align*} 
b< \lambda_i(A_n), \qquad  & \text{for} \quad  i=1,...,q\\
a\leq \lambda_j(A_n)\leq b, \qquad  & \text{for} \quad  j=q+1,...,n,
\end{align*}
then $k$ iterations of CG are sufficient in order to obtain a relative error of $\epsilon$ of the solution, for any initial guess, i.e.
$$ \|\mathbf{e}^k\|_{A_n}/\|\mathbf{e}^0\|_{A_n}\leq\epsilon,$$
where $\mathbf{e}^k$ is the error vector at the $k$th iteration and 
$$ k = q + \lceil \log(2\epsilon^{-1})/\log(\alpha^{-1}) \rceil,$$
with $\alpha=\frac{\sqrt{b}-\sqrt{a}}{\sqrt{b}+\sqrt{a}}$.
\end{lemma}
This lemma can be useful if $q\ll n$, i.e. if there are $q$ outliers in the spectrum of $A_n$, all larger than $b$. Interestingly, the convergence of the CG method does not depend on the magnitude of the outliers, but only on $q$. 
The latter is an explanation of why the convergence rate of the CG applied to a linear system with matrix $\tau^{-1}A_n$ does not significantly depend on $\tau$, as we will show in the numerical section. Indeed, the matrix $\tau^{-1}A_n$ can be written as 
\[
    \tau^{-1}A_n = B_n + \tilde{R}_n
\]
with
\begin{equation}\label{eq:A_eq_B_plus_R}
    B_n = 
    \begin{bmatrix}
	0 & A_{e}^{\mathrm{in},\Gamma} & 0 & 0  \\
	A_{e}^{\Gamma,\mathrm{in}} & A_{e}^{\Gamma} & 0 & 0  \\
        0 & 0 & 0 & A_{i}^{\mathrm{in},\Gamma}  \\
        0 & 0 & A_{i}^{\Gamma,\mathrm{in}} & A_{i}^{\Gamma}  \\
    \end{bmatrix},
    \qquad \tilde{R}_n = \tau^{-1}
    \begin{bmatrix}
	0 & 0 & 0 & 0  \\
	0 &  M_e^\Gamma & 0 & T_{ei}^\Gamma  \\
        0 & 0 & 0 & 0  \\
        0 & (T_{ei}^\Gamma)^T & 0 & M_i^\Gamma  \\
    \end{bmatrix},
\end{equation}
where the first term does not depend on $\tau$ and the second term has rank at most $2N_\Gamma$. Thanks to the Cauchy interlacing theorem (see \cite{MR1477662}), if we choose $[a,b]=[\lambda_{\min}(B_n),\lambda_{\max}(B_n)]$ in Lemma \ref{lemma}, then $q \leq 2N_\Gamma$. In practice, for $\tau$ sufficiently small, we have $q = 2N_\Gamma$, as we show in the numerical section (cf. Figure~\ref{fig:eigs}). Nevertheless, due to the proximity of the lower bound $a$ to zero, the convergence of the CG method is far from satisfactory, and 
it becomes essential to design an appropriate preconditioner for $A_n$.


\subsubsection*{Preconditioning}
We propose a theoretical preconditioner, which is proven to be effective through the spectral analysis carried out in the previous sections. The process of implementing a practical preconditioning strategy will be discussed at the end of the current section.

Denoting by $I_{e}^{\Gamma}$ and $I_{i}^{\Gamma}$ the identity matrices of the same size as $A_{e}^{\Gamma}$ and $A_{i}^{\Gamma}$ respectively, a first observation is that the block diagonal matrix
\[
P_n=
\begin{bmatrix}
        \tau A_{e}^{\mathrm{in}} & 0 & 0 & 0  \\
	0 & I_{e}^{\Gamma} & 0 & 0  \\
        0 & 0 & \tau A_{i}^{\mathrm{in}} & 0  \\
        0 & 0 & 0 & I_{i}^{\Gamma} \\
\end{bmatrix}
\]
is such that the Hermitian matrix
\[
A_n-P_n=
    \begin{bmatrix}
	0 & \tau A_{e}^{\mathrm{in},\Gamma} & 0 & 0  \\
	\tau A_{e}^{\Gamma,\mathrm{in}} & \tau A_{e}^{\Gamma} + M_e^\Gamma - I_{e}^{\Gamma} & 0 & T_{ei}^\Gamma  \\
        0 & 0 & 0 & \tau A_{i}^{\mathrm{in},\Gamma}  \\
        0 & (T_{ei}^\Gamma)^T & \tau A_{i}^{\Gamma,\mathrm{in}} & \tau A_{i}^{\Gamma} + M_i^\Gamma - I_{i}^{\Gamma}  \\
    \end{bmatrix}
\]
has rank less than or equal to $4N_\Gamma=o(n)$, $n=N_i+N_e$, with a similar argument of Theorem \ref{th:general} proof.
In fact the typical behavior is $N_\Gamma\sim \sqrt{n}$ if $\Omega$ is a two-dimensional domain. In the general setting of a $d$-dimensional domain and with any reasonable discretization scheme we have  $N_\Gamma\sim n^{d-1\over d}$, $d\ge 1$.

The function $f^{\rm{square}}$ in formula (\ref{eq:f_square}) is nonnegative, implying that the matrices $T_n(f^{\rm{square}})$ are Hermitian positive definite according to \cite[Theorem 3.1]{GS-II}. The restriction operators $\Pi_{\Omega_i}$ and $\Pi_{\Omega_i}^T$ have full rank and so $\Pi_{\Omega_i}T_n(f^{\rm{square}})\Pi_{\Omega_i}^T$ is Hermitian positive definite. Moreover, in our case the boundary conditions do not alter the positive definiteness and the same holds for $A_{e}^{\mathrm{in}}$.
Hence, the matrix $P_n^{-1/2}A_nP_n^{-1/2}$ is well-defined and we can write
\[
    P_n^{-1/2}A_nP_n^{-1/2} = I_n-P_n^{-1/2}(P_n-A_n)P_n^{-1/2}
\]
with $P_n^{-1/2}(P_n-A_n)P_n^{-1/2}$ having rank less than or equal to $4N_\Gamma=o(n)$, $n=N_i+N_e$, which means that the sequence $\{P_n^{-1/2}(P_n-A_n)P_n^{-1/2}\}_n$ is zero-distributed. Since all the involved matrices are Hermitian,  Exercise 5.3 in book \cite{GS-I} tells us that the matrix sequence $\{P_n^{-1/2}A_nP_n^{-1/2}\}_n$ is distributed as 1 in the eigenvalue sense. Therefore ${P_n}$ could serve as an effective preconditioner for ${A_n}$, also because, as already pointed out, the number of possible outliers cannot grow more than $O(n^{d-1\over d})$, $d\ge 1$, $d$ dimensionality of the physical domain.

We now face the task of finding an efficient method for approximating the matrix-vector products with the inverses of $A_{e}^{\mathrm{in}}$ and $A_{i}^{\mathrm{in}}$. Given the multilevel structure of these matrices and their role as discretization of the Laplacian operator, multigrid techniques are a natural choice. Defining $k=\lfloor n/4 \rfloor$, when dealing with a bi-level Toeplitz matrix of size $n$ generated by $f^{\rm{square}}$, the conventional bisection and linear interpolation operators are an effective choice for restriction and prolongation operators. These can be paired with the Jacobi smoother, which requires suitably chosen relaxation parameters. See \cite{MR2114334} for the theoretical background. A similar multigrid strategy can be applied to the ``subdomain'' matrices $A_{e}^{\mathrm{in}}$ and $A_{i}^{\mathrm{in}}$ thanks to \cite{MR3210157}.

In practice, we will apply one multigrid iteration as preconditioner for the whole matrix $A_n$, avoiding the construction of $P_n$. It has been proven \cite{MR3210157} that if two linear systems have spectrally equivalent coefficient matrices if a multigrid method is effective for a linear system, then it is effective also for the other. The matrix sequences $\{A_n\}_n$ and $\{P_n\}_n$ share the same spectral distribution, but the matrices may not be spectrally equivalent due to outliers, which are at most $4N_\Gamma$ by the Cauchy interlacing theorem. However, outliers are controlled by the CG method, as shown by Lemma  \ref{lemma}.
For a recent reference on spectrally equivalent matrix sequences and their use to design fast iterative solvers see \cite{Ax-last}.

\section{Numerical results}\label{sec:num}

\subsection{Problem and discretization settings}

\begin{figure}
    \centering
    \includegraphics[width=0.4\textwidth]{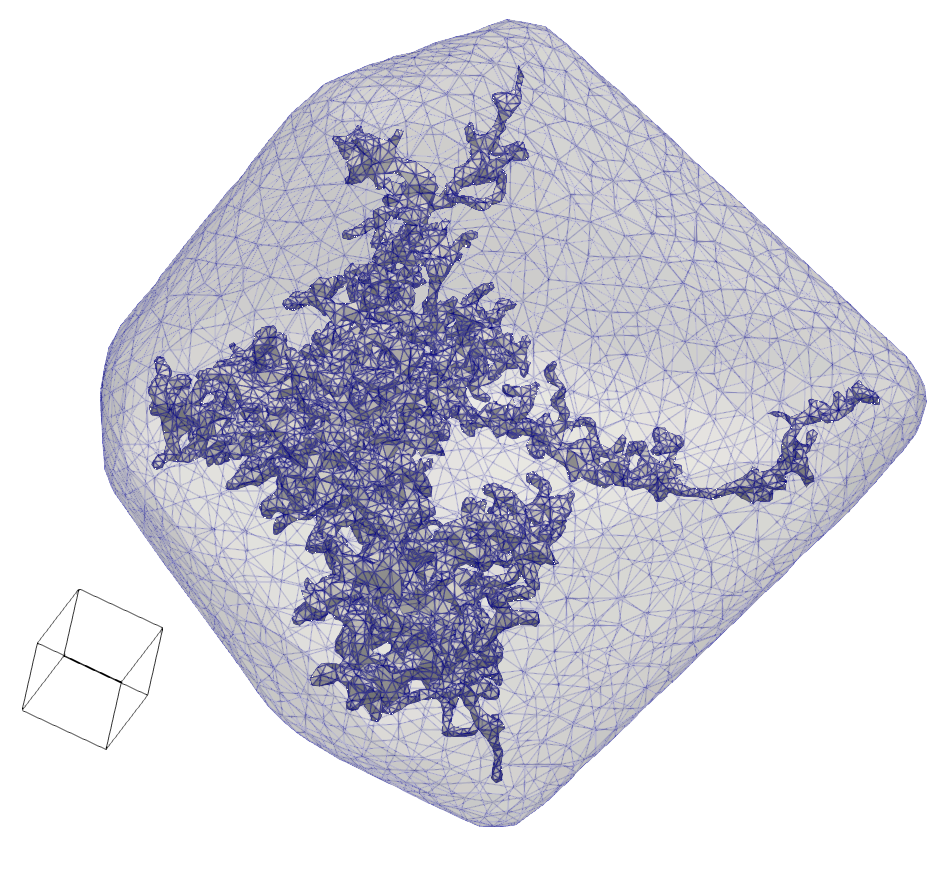} \quad
    \includegraphics[width=0.5\textwidth]{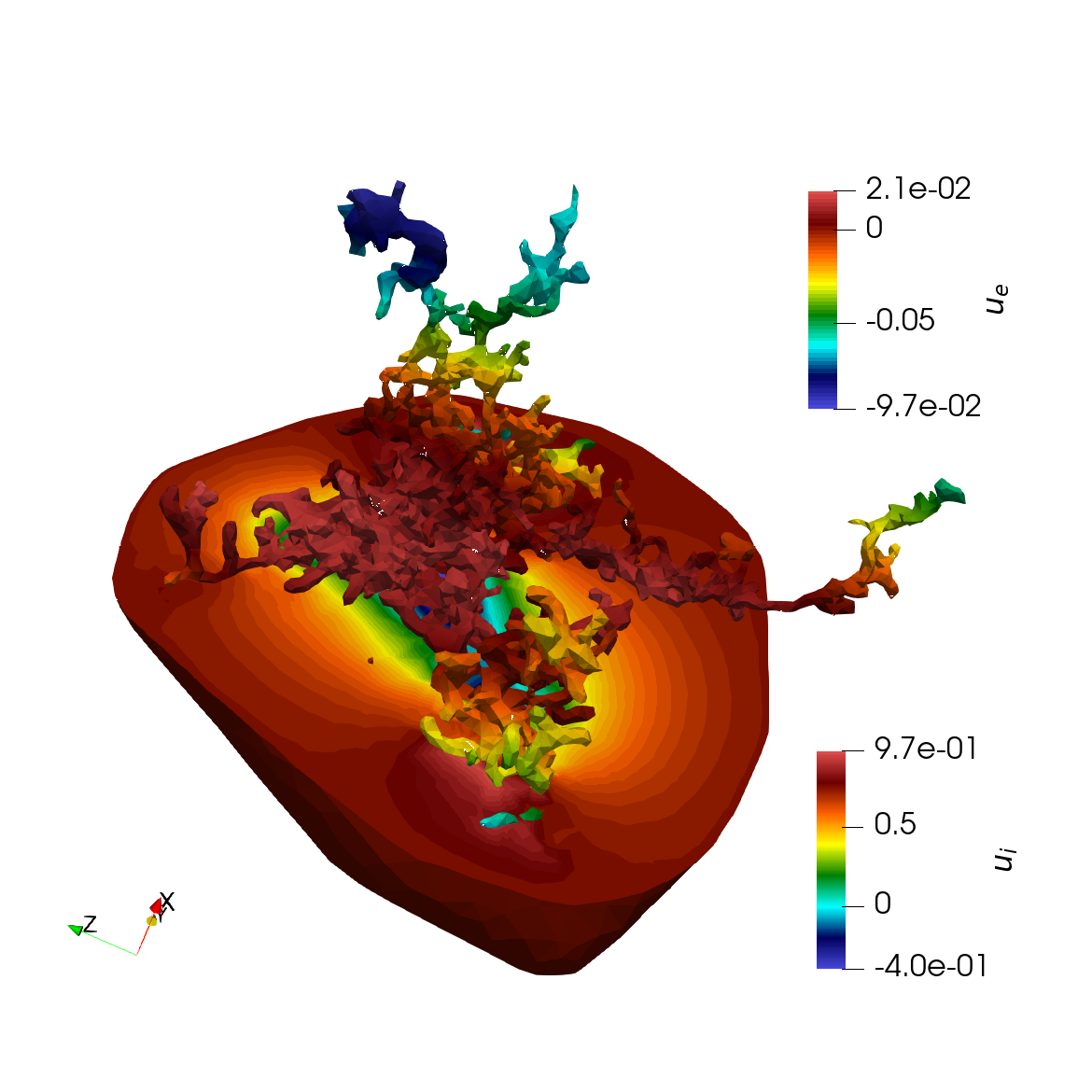}
    \caption{Left: astrocyte geometry and corresponding mesh. A cube with 0.1 sides is shown as a reference. Right: solution corresponding to Table~\ref{tab:astro_ii}. }
    \label{fig:astro}
\end{figure}

We consider the following two scenarios: 

\begin{itemize}
    
    \item[(i)] An idealized geometry for a single cell, i.e. a two-dimensional square domain $\Omega=[0,1]^2$ with $\Omega_i=[0.25,0.75]^2$, discretized with a uniform mesh, cf. Figure~\ref{fig:EMI}. The domain $\Omega$ is discretized considering a uniform grid with $N\times N$ elements. This tessellation results in $n=N_i+N_e=(Np+1)^2+2Np$ degrees of freedom\footnote{The term $2Np=N_\Gamma$ is present since the membrane degrees of freedom are repeated both for $\Omega_i$ and $\Omega_e$.}, with $p$ the finite element order.

    \item[(ii)] A three-dimensional geometry, representing one astrocyte \cite{abdellah2023ultraliser}, cf. Figure~\ref{fig:astro} discretized with 32365 grid nodes leading to $n=212548$ degrees of freedom for $p = 2$ (as used in the numerical experiments).

\end{itemize}
In both cases, we set $\sigma_i=\sigma_e=1$, pure Dirichlet boundary conditions, i.e. $\partial\Omega_D=\partial\Omega$ in \eqref{bc1}, and the following right-hand side source in \eqref{eq::f}:
   $$ f(x,y) = \sin(2\pi x) \sin(2\pi y).$$ 
%

\subsection{Implementation and solver parameters }
We use FEniCS \cite{alnaes2015fenics,logg2012automated} for parallel finite element assembly and \texttt{multiphenics}\footnote{\url{https://multiphenics.github.io/index.html}} to handle multiple meshes with common interfaces and the corresponding mixed-dimensional weak forms. FEniCS wraps PETSc parallel solution strategies. For multilevel preconditioning, we use \texttt{hypre} algebraic multigrid (AMG) \cite{falgout2002hypre}, with default options, cf. Appendix~\ref{appendix_ksp}. 
For comparative studies, we use MUMPS as a direct solver and incomplete LU (ILU) preconditioning with zero fill-in.
Parallel ILU relies on a block Jacobi preconditioned, approximating each block inverse with ILU. For iterative strategies, we use a tolerance of $10^{-6}$ for the relative unpreconditioned residual, as stopping criteria.

The implementation is verified considering the same benchmark problem as in \cite[Section 3.2]{tveito2017cell}, obtaining the expected convergence rates.

\subsection{Experiments: eigenvalues and CG convergence}
\begin{figure}
    \centering
    \includegraphics[width=\textwidth]{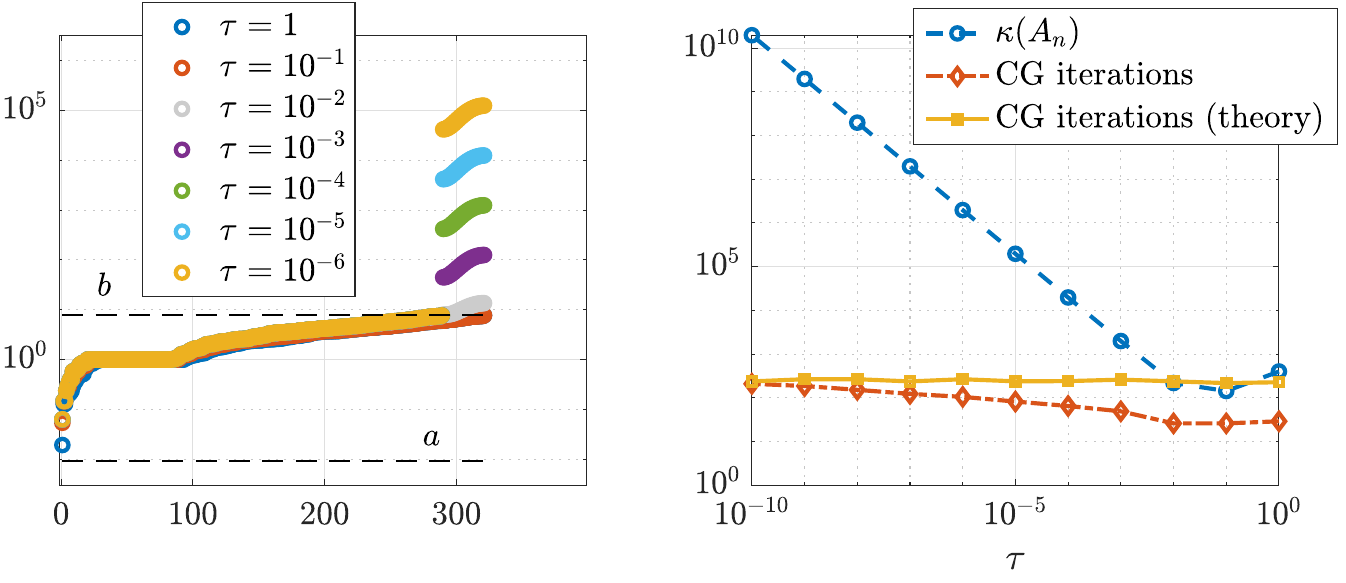}
    \caption{Left: spectrum of $\tau^{-1}A_n$ varying $\tau$ for $(N,p)=(16,1)$. The $[a,b]=[\lambda_{\min}(B_n), \lambda_{\max}(B_n)]$ interval is shown, with reference to Lemma~\ref{lemma} and equation \eqref{eq:A_eq_B_plus_R}. Crucially, the majority of eigenvalues are included in $[a,b]$, which is independent of $\tau$. Right: condition number of $A_n$ and CG iterations to convergence. In yellow the theoretical bound is shown, which depends on $a,b$ and the relative error $\epsilon$.}
    \label{fig:eigs}
\end{figure}

We consider the behaviour of unpreconditioned CG for $\tau \to 0$ and scenario (i). For simplicity, in this section, we use a unit right-hand side, i.e. $\ff=\mathbf{1}/\|\mathbf{1}\|_2$ in \eqref{eq:system}.
In Figure~\ref{fig:eigs} an example of the spectrum of $A_n$ is shown; as presented in Section~\ref{ssec:preconditioner}, varying $\tau$ influences only $2N_\Gamma$ eigenvalues. Using Lemma~\ref{lemma}, with $q = 2N_\Gamma$, results in a strict upper bound for the number of CG iterations, which is observed in practice, despite the increasing condition number $\kappa(A_n)$.

\begin{figure}
    \centering    
    \includegraphics[width=0.48\textwidth]{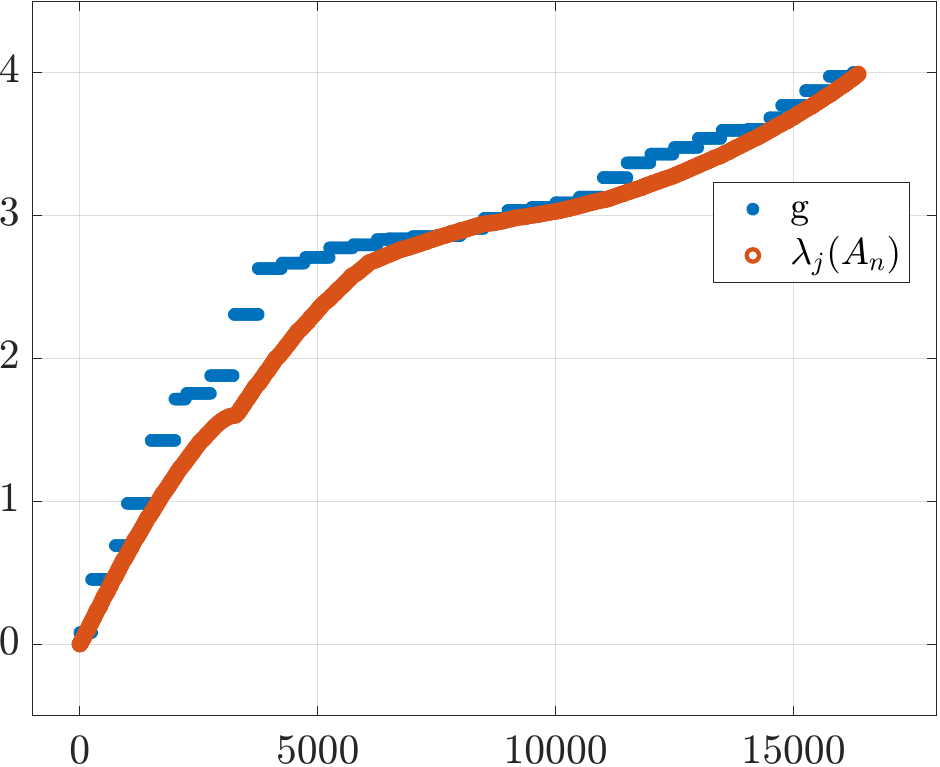}
    \includegraphics[width=0.48\textwidth]{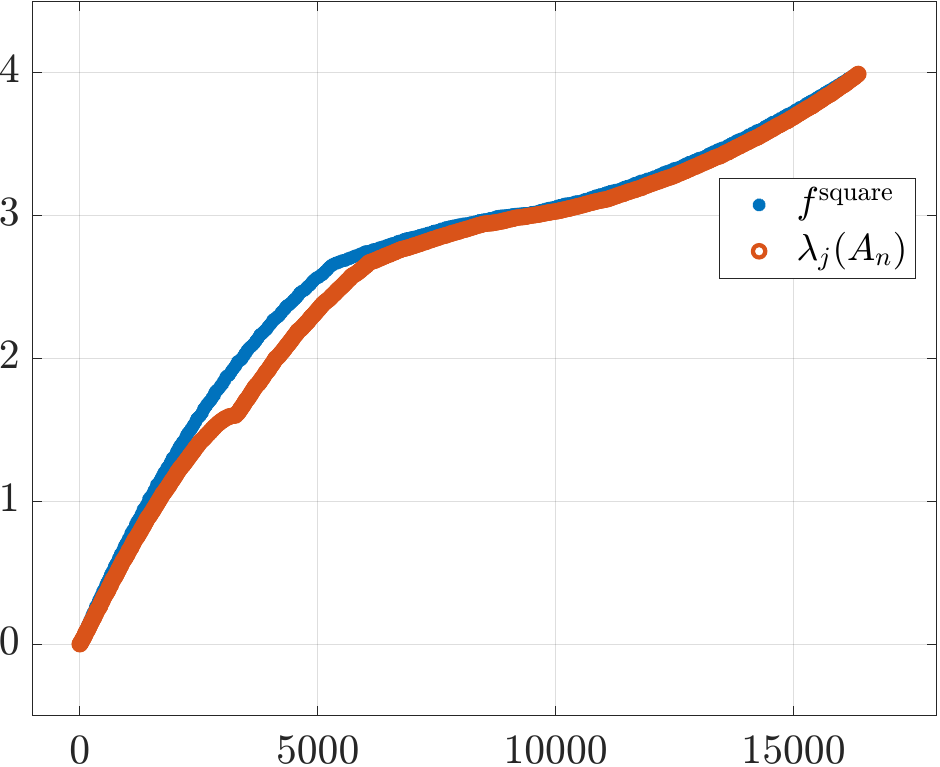}
    \caption{Left: comparison between the eigenvalues $\lambda_j(A_{n})$ and the samples of $g$ over a uniform grid on $D$.
    Right: comparison between the eigenvalues  $\lambda_j(A_{n})$ and the samples of $f^{\rm{square}}$ over a uniform grid on $[0,\pi]^2$. In both cases $(N,p)=(128,1)$ and $\tau= 1$. 
    }
    \label{fig:sampling_g}
\end{figure}

In this idealized setting, imposing $\tau = 1$, we numerically confirm the results of sections~\ref{ssec:symbol_analysis}-\ref{ssec:stiffness_matrices}, showing that the spectral distribution
\[
    \{A_n\}_n\sim_\lambda (g,D), 
\]
with the symbol domain  $$D=[0,1] \times [0,\pi]^2 \times [0,\pi]^2$$ 
and
    $$g(x,\boldsymbol{\theta}_i,\boldsymbol{\theta}_e)=f^{\rm{square}}(\boldsymbol{\theta}_i) \psi_{[0,\frac{1}{4}]}(x) + f^{\rm{square}}(\boldsymbol{\theta}_e) \psi_{(\frac{1}{4},1]}(x),$$
reasonably approximates the eigenvalues of $A_n$, given a uniform sampling even for very few grid points in each of the 5 dimensions of $D$. In order to make visualization possible, we evaluate the function $g$ on its domain and then we arrange the evaluations in ascending order, which corresponds to considering the monotone rearrangement of $g$. Note that in this setting the number $r$ defined in Theorem \ref{th:general} is equal to $1/4$. 
In the left panel in Figure~\ref{fig:sampling_g} we observe that the spectrum of $A_{n}$ is qualitatively described by the samplings of the function $g$ over a uniform grid on $D$. Moreover, in our particular setting, also the spectral distribution 
\[
    \{A_n\}_n\sim_\lambda (f^{\rm{square}},[0,\pi]^2)
\]
holds, and we numerically confirm this result in the right panel of Figure~\ref{fig:sampling_g}, where the eigenvalues of $A_{n}$ are described by the samplings of the function $f^{\rm{square}}$ over a uniform grid on $[0,\pi]^2$.

\subsection{Experiments: preconditioning}
In Table~\ref{tab:CG} we report CG iterations to convergence for scenario (i), varying the discretization parameters $N,p$ and $\tau$. 
Similarly, in Table~\ref{tab:PCG}, we report iterations to convergence using the AMG preconditioned CG (PCG), showing robustness w.r.t. all discretization parameters. We remark that, for PCG, time-to-solution depends linearly on $n$, as expected for multilevel methods.
In Table~\ref{tab:astro_i} we report runtimes and iterations to convergence for various preconditioning strategies for scenario (i); assembly and direct solution timings are also reported for practicality. In terms of efficiency, the AMG and ILU preconditioners are preferable. In Table~\ref{tab:astro_ii}, we show strong scaling data corresponding to scenario (ii), i.e. the astrocyte geometry. In this case, the ILU preconditioner is preferable both in terms of runtime and parallel efficiency.

Summing up the numerical results, we can observe how the AMG preconditioned CG exhibit a near-to-optimal robustness, converging in 5 to 7 iterations for all the cases considered. For time-to-solution, an ILU preconditioning might be preferable for three-dimensional complex geometries, while AMG remains the most convenient choice for structured discrete problems.  Let us remark that the ILU preconditioner can be particularly suitable for multiple time steps since the ILU decomposition has to be computed only once.

\begin{table}[]
    \centering
    \begin{tabular}{l|c|c|c|c|c}
         $(N,p)$ & $(32,1)$ & $(64,1)$ & $(128,1)$ & $(256,1)$ & $(512,1)$\\
         \hline
         $n$ & 1153 & 4353 & 16897 & 66561 & 264193\\
         \hline                  
         $\tau=1$       & 49 & 95 & 186 & 365 & 718\\
         $\tau=10^{-1}$ & 45 & 90 & 176 & 346 & 686 \\
         $\tau=10^{-2}$ & 41 & 79 & 151 & 296 & 583\\
         $\tau=10^{-3}$ & 71 & 122 & 185 & 263 & 487\\
          \\
            \\
         $(N,p)$ & $(16,2)$ & $(32,2)$ & $(64,2)$ & $(128,2)$ & $(256,2)$ \\
         \hline        
         $\tau=1$       & 56 & 110 & 215 & 422 & 830 \\
         $\tau=10^{-1}$ & 50 & 104 & 204 & 400 & 793 \\
         $\tau=10^{-2}$ & 47 & 92  & 175 & 343 & 674\\
         $\tau=10^{-3}$ & 79 & 129 & 192 & 298 & 563
    \end{tabular}
    \caption{Number of CG iterations to convergence, scenario (i).}
    \label{tab:CG}
\end{table}

\begin{table}[]
    \centering
    \begin{tabular}{l|c|c|c|c|c}
         $(N,p)$ & $(32,1)$ & $(64,1)$ & $(128,1)$ & $(256,1)$ & $(512,1)$\\
         \hline
         $n$ & 1153 & 4353 & 16897 & 66561 & 264193\\
         \hline                  
         $\tau=1$       & 5 & 5 & 5 & 6 & 5 \\
         $\tau=10^{-1}$ & 5 & 5 & 6 & 6 & 6 \\
         $\tau=10^{-2}$ & 5 & 5 & 5 & 6 & 7 \\
         $\tau=10^{-3}$ & 5 & 5 & 5 & 6 & 6 \\
          \\
            \\
         $(N,p)$ & $(16,2)$ & $(32,2)$ & $(64,2)$ & $(128,2)$ & $(256,2)$ \\
         \hline        
         $\tau=1$       & 5 & 5 & 5 & 6 & 5 \\
         $\tau=10^{-1}$ & 5 & 6 & 5 & 6 & 6 \\
         $\tau=10^{-2}$ & 4 & 5 & 5 & 5 & 6 \\
         $\tau=10^{-3}$ & 6 & 5 & 5 & 6 & 6 
    \end{tabular}
    \caption{Number of PCG iterations to convergence, scenario (i), using a single AMG iteration as preconditioner.}
    \label{tab:PCG}
\end{table}

\begin{table}[]
    \centering
      \begin{tabular}{l|c}
          & Runtime [its.] \\
         \hline
         Assembly & 2.6  \\
         Direct   & 1.6 \\
         CG       & 1.4 [583] \\         
         PCG(Jacobi)  & 1.5 [557] \\         
         PCG(SOR)     & 1.2 [199] \\               
         PCG(ILU) & 0.7 [120] \\ 
         PCG(AMG)     & 0.3 [5]             
    \end{tabular}
    \caption{Runtimes (s) and iterations to convergence (in square brackets) for various solution strategies for $\tau = 0.01$ and scenario (i) with $(N,p)=(512,1)$ for a single core. For completeness, finite element assembly timing is also reported.}
    \label{tab:astro_i}
\end{table}

\begin{table}[]
    \centering
        
    \begin{tabular}{l|c|c|c|c}
         Cores & 1 & 2 & 4 & 8\\
         \hline
         Assembly & 8.4 & 5.5 & 3.4 & 2.3 \\
         Direct   & 16 & 12.1 & 10.6 & 9.9 \\
         CG       & 6.3 [725] & 3.5 [724] & 2.2 [724] & 2.0 [720] \\         
         PCG(Jacobi)  & 3.0 [337] & 1.7 [337] & 1.0 [337]  & 1.0 [337]  \\
         PCG(SOR)   & 2.6  [129] & 2.0 [174] & 1.4  [181] & 1.2 [185] \\         
         PCG(ILU) & 2.2 [98] & 1.3 [104] & 0.8 [113] & 0.6 [119]  \\  
         PCG(AMG)      & 5.9 [7] & 3.9 [7] & 3.0 [7]      & 2.9 [7]                   
      \end{tabular}
    \caption{Runtimes (s) and iterations to convergence (in square brackets) for various parallel solution strategies for $\tau = 0.01$ and scenario (ii), cf. the astrocyte geometry in Figure~\ref{fig:astro}, and $p = 2$. For completeness, finite element assembly timings are also reported.}
    \label{tab:astro_ii}
\end{table}

\section{Concluding remarks}\label{sec:final}

We described numerical approximation schemes for the EMI equations and studied the structure and spectral features of the coefficient matrices obtained from a finite element discretization in the so-called \textit{single-dimensional} formulation. The obtained spectral information has been employed for designing appropriate (preconditioned) Krylov solvers. Numerical experiments have been presented and critically discussed; depending on the context (dimensionality, geometry) the CG method, preconditioned with AMG or ILU results in an efficient, scalable, and robust solution strategy. The spectrum analysis made it possible to understand the convergence properties, somehow counterintuitive, of CG in this context.

Due to the intrinsic generality of the tools in Sections~\ref{ssec:GLTbackground}-\ref{ssec:symbol_analysis}, the given machinery can be applied to a variety of discretization schemes including isogeometric analysis, finite differences, and finite volumes, in the spirit of the sections of books \cite{GS-II,block-glt-dD} dedicated to applications.


\subsection*{Acknowledgments} Pietro Benedusi and Marie Rognes acknowledge support from the Research Council of Norway via FRIPRO grant \#324239 (EMIx) and from the national infrastructure for computational science in Norway, Sigma2, via grant \#NN8049K. 
 Paola Ferrari and Stefano Serra-Capizzano are partially supported by the Italian Agency INdAM-GNCS. Furthermore, the work of Stefano Serra-Capizzano is funded from the European High-Performance Computing Joint Undertaking  (JU) under grant agreement No 955701. The JU receives support from the European Union’s Horizon 2020 research and innovation programme and Belgium, France, Germany, and Switzerland.
 Stefano Serra-Capizzano is also grateful for the support of the Laboratory of Theory, Economics and Systems – Department of Computer Science at Athens University of Economics and Business. Special thanks are extended to J\o rgen Dokken, Miroslav Kuchta, Francesco Ballarin, Lars Magnus Valnes, Halvor Herlyng, Abdellah Marwan, and Marius Causemann for their precious help.

\bibliographystyle{plain}
\bibliography{references}

\appendix
\section{Appendix: multigrid solver parameters}\label{appendix_ksp}
We provide the default AMG parameters in PETSc obtained using the \texttt{ksp\_view} command.

\vspace{10pt}
\noindent
{\tt KSP Object: 1 MPI process \\
\null\quad type: cg \\
\null\quad maximum iterations=10000, initial guess is zero\\
 \null\quad tolerances:  relative=1e-06, absolute=1e-50, \null\quad divergence=10000.\\
  \null\quad left preconditioning\\
  \null\quad using UNPRECONDITIONED norm type for convergence test\\
PC Object: 1 MPI process\\
 \null\quad type: hypre\\
  \null\qquad  HYPRE BoomerAMG preconditioning\\
   \null\qquad\quad   Cycle type V\\
    \null\qquad\quad  Maximum number of levels 25\\
      \null\qquad\quad Maximum number of iterations PER hypre call 1\\
      \null\qquad\quad Convergence tolerance PER hypre call 0.\\
      \null\qquad\quad Threshold for strong coupling 0.25\\
      \null\qquad\quad Interpolation truncation factor 0.\\
      \null\qquad\quad Interpolation: max elements per row 0\\
      \null\qquad\quad Number of levels of aggressive coarsening 0\\
      \null\qquad\quad Number of paths for aggressive coarsening 1\\
      \null\qquad\quad Maximum row sums 0.9\\
      \null\qquad\quad Sweeps down         1\\
      \null\qquad\quad Sweeps up           1\\
      \null\qquad\quad Sweeps on coarse    1\\
      \null\qquad\quad Relax down          symmetric-SOR/Jacobi\\
      \null\qquad\quad Relax up            symmetric-SOR/Jacobi\\
      \null\qquad\quad Relax on coarse     Gaussian-elimination\\
      \null\qquad\quad Relax weight  (all)      1.\\
      \null\qquad\quad Outer relax weight (all) 1.\\
      \null\qquad\quad Using CF-relaxation\\
      \null\qquad\quad Not using more complex smoothers.\\
      \null\qquad\quad Measure type        local\\
      \null\qquad\quad Coarsen type        Falgout\\
      \null\qquad\quad Interpolation type  classical\\
      \null\qquad\quad SpGEMM type         cusparse\\
  linear system matrix = precond matrix\\
}

\end{document}